\documentclass[10pt]{amsart}
\usepackage{amsmath, amsthm, amssymb}

\theoremstyle{plain}
\newtheorem{thm}{Theorem}

\title{More Counterexamples to the Unit Conjecture for Group Rings}
\author{Alan G. Murray}
\address{IDA Center for Computing Sciences, Bowie, MD}
\email{murray@super.org}

\begin{document}

\begin{abstract}
Extending the discovery by Giles Gardam of a concrete counterexample
to Kaplansky's unit conjecture in characteristic 2, a family of
counterexamples for every prime characteristic is presented.
\end{abstract}

\maketitle

\section{Background}

Consider the group ring $K[G]$, where $K$ is a field and $G$ is a
torsion-free group. The unit conjecture claims that the only units of
$K[G]$ are non-zero scalar multiples of elements of $G$. The conjecture
was posed in 1970 by Kaplansky, and is commonly attributed to him,
although it was first stated in a 1940 thesis by Higman.
Then in February 2021, Giles Gardam concluded an online talk by presenting
a counterexample for $K=\mathbb{F}_2$, disproving the conjecture and
opening a new direction of research.

This note follows this direction established by Gardam \cite{Gar21}, using much
of the same notation but with a few variations. Let $P$ be the torsion-free group
$\langle a,b~|~(a^2)^b = a^{-2}, (b^2)^a = b^{-2}\rangle$ and set
$x=a^2$, $y=b^2$, $z=(ab)^2$. When dealing with functions
$f(x,y,z)$, it is useful to denote inversions of variables like this:
\begin{displaymath}
f_x(x,y,z)=f(x^{-1},y,z), \quad f_{xy}(x,y,z)=f(x^{-1},y^{-1},z), \quad {\rm etc.}
\end{displaymath}
Using this notation, we have the relations
$$af=f_{yz}a, \quad bf=f_{xz}b, \quad abf=f_{xy}ab.$$
Combined with the identity
$ba=x^{-1}yz^{-1}ab$, we can write any element of $P$ uniquely in the form
$f(x,y,z)g$, where $g\in \{1,a,b,ab\}$.

The condition $(p'+q'a+r'b+s'ab)(p+qa+rb+sab)=1$ is equivalent to
\begin{displaymath}
\begin{array}{rcrcrcrc}
p'p&+&xq'q_{yz}&+&yr'r_{xz}&+&zs's_{xy}& = 1 \\
p'q&+&q'p_{yz}&+&x^{-1}z^{-1}r's_{xz}&+&y^{-1}s'r_{xy}& = 0 \\
p'r&+&xq's_{yz}&+&r'p_{xz}&+&y^{-1}zs'q_{xy}& = 0 \\
p's&+&q'r_{yz}&+&x^{-1}yz^{-1}r'q_{xz}&+&s'p_{xy}& = 0
\end{array}
\end{displaymath}
Using Gardam's selection
$(p',q',r',s')=(x^{-1}p_{yz},-x^{-1}q,-y^{-1}r,z^{-1}s_{yz})$, 
\begin{displaymath}
\begin{array}{rcrcrcrc}
x^{-1}p_{yz}p&+&-qq_{yz}&+&-rr_{xz}&+&s_{yz}s_{xy}& = 1 \\
&&&&-x^{-1}z^{-1}y^{-1}rs_{xz}&+&y^{-1}z^{-1}s_{yz}r_{xy}& = 0 \\
x^{-1}p_{yz}r&+&-qs_{yz}&+&-y^{-1}rp_{xz}&+&y^{-1}s_{yz}q_{xy}& = 0 \\
x^{-1}p_{yz}s&+&-x^{-1}qr_{yz}&+&-x^{-1}z^{-1}rq_{xz}&+&z^{-1}s_{yz}p_{xy}& = 0
\end{array}
\end{displaymath}
After choosing $p_{xy}=x^{-1}y^{-1}p$, $q_{xy}=yq$,
$r_{xy}=x^{-1}r$, $s_{xy}=s$, the second and third equations
are satisfied, leaving
\begin{displaymath}
\begin{array}{cccccccc}
  pp_{yz}&+&xss_{yz}&=&rr_{yz}&+&xqq_{yz}&+~x \\
  ps_{yz}&+&yzp_{yz}s&=&rq_{yz}&+&yzqr_{yz}&
\end{array}
\end{displaymath}
From here, it is relatively easy to verify Gardam's result
for characteristic 2.
\begin{thm}[Gardam]
Let $P$ be the torsion-free group
$\langle \, a, b \, | \, (a^2)^b=a^{-2}, \, (b^2)^a=b^{-2} \, \rangle$
and set $x=a^2$, $y=b^2$, $z=(ab)^2$. Set
\begin{displaymath}
\begin{array}{rl}
p &= (1+x)(1+y)(1+z^{-1}) \\
q &= x^{-1}y^{-1}+x+y^{-1}z+z \\
r &= 1+x+y^{-1}z+xyz \\
s &= 1 + (x+x^{-1}+y+y^{-1})z^{-1}.
\end{array}
\end{displaymath}
Then $p+qa+rb+sab$ is a non-trivial unit in the group ring $\mathbb{F}_2[P]$.
\end{thm}

To find additional solutions, it suffices to match the two left
hand sides (functions of $p$ and $s$) and the two right hand sides
(functions of $r$ and $q$). A simple search over low-degree candidates
(degree at most 1 in $\{x^{\pm 1},y^{\pm 1},z^{\pm 1}\}$) yielded 18 solutions
for $\mathbb{F}_2[P]$, all related to the counterexample
just mentioned. In particular, any of three choices for $(p,s)$
can be paired with any of six choices for $(q,r)$:
\begin{displaymath}
\begin{array}{cl}
(p,s) : & (p,s), (p_z,zs), (p_z,z^{-1}s_z) \\
(q,r) : & (q,r), (q,zr_z), (q_x,r_y), (q_x,zr_{yz}), (z^{-1}q,z^{-1}r), (z^{-1}q_x,z^{-1}r_y).
\end{array}
\end{displaymath}

\section{New Results}

In order to find non-trivial units for the group $P$ over fields of
other characteristic, it is natural to start with characteristic 3.
Repeating the search over low-degree candidates in characteristic 3
yielded two related counterexamples:

\begin{thm}
Let $P$ be the torsion-free group
$\langle \, a, b \, | \, (a^2)^b=a^{-2}, \, (b^2)^a=b^{-2} \, \rangle$
and set $x=a^2$, $y=b^2$, $z=(ab)^2$. Set
\begin{displaymath}
\begin{array}{rl}
p &= (1+x)(1+y)(z^{-1}-z) \\
q &= (1+x)(x^{-1}+y^{-1})(1-z^{-1})+(1+y^{-1})(z-z^{-1}) \\
r &= (1+y^{-1})(x+y)(z-1)+(1+x)(z-z^{-1}) \\
s &= -z + (1+x+x^{-1}+y+y^{-1})(z^{-1}-1).
\end{array}
\end{displaymath}
Then both $p+qa+rb+sab$ and $p+q_xa+r_yb+sab$ are non-trivial units
in the group ring $\mathbb{F}_3[P]$.
\end{thm}

There are evident similarities between the $\mathbb{F}_2[P]$ and
$\mathbb{F}_3[P]$ counterexamples, and this suggests the possibility
of a search over polynomials in $z$ only. So we consider a
counterexample of the following form:
\begin{displaymath}
\begin{array}{rl}
p &= (1+x)(1+y)f_1(z) \\
q &= (1+x)(x^{-1}+y^{-1})f_2(z)+(1+y^{-1})f_3(z) \\
r &= (1+y^{-1})(x+y)f_4(z)+(1+x)f_5(z) \\
s &= (x+2+x^{-1}+y+2+y^{-1})f_6(z)+f_7(z)
\end{array}
\end{displaymath}

Notice that these choices still satisfy $p_{xy}=x^{-1}y^{-1}p$,
$q_{xy}=yq$, $r_{xy}=x^{-1}r$, $s_{xy}=s$. The seven Laurent polynomials
$f_1(z)$ through $f_7(z)$ will henceforth be assumed to be functions
only of $z$, and the argument will be suppressed.
Instead we write $f_1=f_1(z)$, $f_1^*=f_1(z^{-1})$, etc. The condition
$(p'+q'a+r'b+s'ab)(p+qa+rb+sab)=1$ is equivalent to
\begin{displaymath}
\begin{array}{rl}
f_7^*f_7 &= 1 \\
f_3^*f_3 &= f_5^*f_5 = (f_7^*f_6+f_6^*f_7) + 4f_6^*f_6 \\
f_2^*f_2 &= f_4^*f_4 = f_6^*f_6 \\
f_1^*f_1 &= f_2^*f_3 + f_4^*f_5 \\
f_2^*f_5 &= zf_5^*f_2 = f_3^*f_4 = zf_4^*f_3 = (f_6^*f_1+zf_1^*f_6)-(f_2^*f_4+zf_4^*f_2) \\
f_3^*f_5+zf_5^*f_3&=(f_7^*f_1+zf_1^*f_7) + 4(f_2^*f_4+zf_4^*f_2)
\end{array}
\end{displaymath}

It can be shown that assuming any of the seven polynomials is zero leads
to either a contradiction or a trivial solution. Thus we may assume all
seven polynomials are nonzero. For the purpose of just finding one
solution, not all solutions, let's assume $f_4=z^tf_2$. Then we can
deduce $f_5=z^{1-t}f_3$. This reduces the system of equations to this:
\begin{displaymath}
\begin{array}{rl}
f_7^*f_7 &= 1 \\
f_3^*f_3 &= (f_7^*f_6+f_6^*f_7) + 4f_6^*f_6 \\
f_2^*f_2 &= f_6^*f_6 \\
f_1^*f_1 &= (1+z^{1-2t})f_2^*f_3 \\
f_2^*f_3 &= z^{2t-1}f_3^*f_2 = (f_6^*f_1+zf_1^*f_6)-(z^t+z^{1-t})f_2^*f_2 \\
(z^t+z^{1-t})f_3^*f_3 &= (f_7^*f_1+zf_1^*f_7) + 4(z^t+z^{1-t})f_2^*f_2
\end{array}
\end{displaymath}
Combining the third, fourth and fifth equations, we can deduce
\begin{displaymath}
\left(f_1-\left(z^t+z^{1-t}\right)f_6\right)\cdot\left(f_1^*-\left(z^{-t}+z^{t-1}\right)f_6^*\right) = 0,
\end{displaymath}
so we know $f_1=(z^t+z^{1-t})f_6$. Plugging this into the fifth equation and
simplifying with the third equation, we get $f_3 = (1+z^{2t-1})f_2$.
Incorporating this, the system reduces to
\begin{displaymath}
\begin{array}{rl}
f_7^*f_7 &= 1 \\
z^{2t-1}\left(1-z^{1-2t}\right)^2\cdot f_6^*f_6 &= f_7^*f_6+f_6^*f_7 \\
f_2^*f_2 &= f_6^*f_6
\end{array}
\end{displaymath}
from which we conclude
\begin{displaymath}
\left(f_7-z^{2t-1}\left(1-z^{1-2t}\right)^2f_6\right)\cdot\left(f_7^*-z^{1-2t}\left(1-z^{2t-1}\right)^2f_6^*\right)=1.
\end{displaymath}
At this point, we are free to set $f_2=z^wf_6$ for any integer $w$,
but in characteristic zero, we still have a problem setting $f_6$ and $f_7$.

However, in characteristic $d>0$, we have the identity
$(1+f(z))^d\equiv 1+f(z)^d$, that we can use to our advantage.
Let $g=z^{1-2t}$ and $f_6=(1-g)^{d-2}$, so that
\begin{displaymath}
f_7-z^{2t-1}(1-g)^2f_6=f_7-z^{2t-1}(1-g)^d\equiv f_7-z^{2t-1}-z^{2t-1}(-g)^d.
\end{displaymath}
Now we can select $f_7=z^{2t-1}$ so that both $f_7$ and
$f_7-z^{2t-1}(1-z)^2f_6$ are monomials, as required. Backsolving,
we get the following choices for the seven Laurent polynomials:
\begin{displaymath}
\begin{array}{rl}
f_1 &= \left(z^t+z^{1-t}\right)\cdot h \\
f_2 &= z^w\cdot h \\
f_3 &= z^w\left(1+z^{2t-1}\right)\cdot h \\
f_4 &= z^{w+t}\cdot h \\
f_5 &= z^w\left(z^t+z^{1-t}\right)\cdot h \\
f_6 &= h \\
f_7 &= z^{2t-1} \\
{\rm where~}~h &= \left(1-z^{1-2t}\right)^{d-2}
\end{array}
\end{displaymath}

Summarizing, we have a family of counterexamples for each prime characteristic:
\begin{thm}
Let $P$ be the torsion-free group
$\langle \, a, b \, | \, (a^2)^b=a^{-2}, \, (b^2)^a=b^{-2} \, \rangle$
and set $x=a^2$, $y=b^2$, $z=(ab)^2$. Choose a prime $d$ and
integers $t$ and $w$, and set 
\begin{displaymath}
\begin{array}{rl}
h &= \left(1-z^{1-2t}\right)^{d-2} \\
p &= (1+x)(1+y)(z^t+z^{1-t})h \\
q &= z^w\left[(1+x)(x^{-1}+y^{-1})+(1+y^{-1})(1+z^{2t-1})\right]h \\
r &= z^w\left[(1+y^{-1})(x+y)z^t+(1+x)(z^t+z^{1-t})\right]h \\
s &= z^{2t-1} + (4+x+x^{-1}+y+y^{-1})h.
\end{array}
\end{displaymath}
Then $p+qa+rb+sab$ is a non-trivial unit in the group ring $\mathbb{F}_d[P]$.
\end{thm}

\end{document}